%% file: NSUM.tex
\documentclass[10pt]{article}

\usepackage{amssymb,amsmath,latexsym}
\voffset
-1.5cm

\usepackage{a4,amscd,graphicx,epsfig}
\usepackage[all]{xy}
\everymath{\displaystyle}

\newtheorem{teo}{Theorem}[section]
\newtheorem{lem}[teo]{Lemma}

\newtheorem{prop}[teo]{Proposition}
\newtheorem{defi}[teo]{Definition}

\newtheorem{conj}[teo]{Conjecture}
\newtheorem{remark}[teo]{Remark}
\newtheorem{remarks}[teo]{Remarks}

\newcommand{\mr}{\mathbb{R}}
\newcommand{\mc}{\mathbb{C}}
\newcommand{\mz}{\mathbb{Z}}
\newcommand{\mh}{\mathbb{H}}

\newcommand{\Dd}{{\mathcal D}}

\newcommand{\Ii}{{\mathcal I}}

\newcommand{\Pp}{{\mathcal P}}

\newcommand{\Tt}{{\mathcal T}}

\newcommand{\cG}{{\mathfrak c}}

\newcommand{\C}{{\mathbb C}}
\newcommand{\Q}{{\mathbb Q}}
\newcommand{\Z}{{\mathbb Z}}
\newcommand{\R}{{\mathbb R}}

\title{QHI Theory, II: \\ Dilogarithmic and Quantum Hyperbolic Invariants \\ 
of 3-Manifolds with $PSL(2,\C)$-Characters\\
{\small [Abridged Version]}}

\author {St\'ephane Baseilhac and Riccardo Benedetti}

\date {}

\begin{document}

\maketitle

\vspace{0.7cm}

\noindent Dipartimento di Matematica, Universit\`a di Pisa, Via
F. Buonarroti, 2, I-56127 PISA. Emails: baseilha@mail.dm.unipi.it, benedett@dm.unipi.it.\\

\noindent The first author benefits from a Marie Curie Fellowship, contract
No MCFI-2001-00605. 

\vspace{1cm}

\input NSUMINTRO.tex

\input SUMIDEAL.tex
\input NSUMDILOG.tex

\input NSUMDINV.tex

\input NSUMQHI.tex
\input SUMSCISSORS.tex
\input SUMCOMMENT.tex

\input NSUMVOL.tex

\input SUMCUSPED.tex

\input SUMBIB.tex
\end{document}

%% file: NSUMINTRO.tex
\begin{abstract}
\noindent We give parallel  constructions of an invariant  ${\rm R}(W, \rho)$, based on the classical
Rogers dilogarithm, and of quantum hyperbolic invariants (QHI), based on the Faddeev-Kashaev
quantum dilogarithms,  for flat $PSL(2,\mc)$-bundles $\rho$ over closed oriented $3$-manifolds $W$. 
All these invariants are explicitely computed as a sum or state sums over the same hyperbolic ideal
tetrahedra of the idealization of any fixed simplicial 1-cocycle description of $(W,\rho)$ of a special kind,
called a $\Dd$-triangulation.  ${\rm R}(W, \rho)$ recovers the volume and the Chern-Simons invariant of $\rho$,
and we conjecture that it determines the semi-classical limit of the QHI.  
\end{abstract}


\section{Introduction} \label{sumintro}

\noindent 
This text is a summary of our recent results. We give here with full
details the complete apparatus of definitions, constructions and
statements, with comments on the key points explaining the logic of
our work. The proofs will be given in a forthcoming expanded version
of this paper. Our aim here is to allow the interested reader to have,
as much as possible, a quick access to these results. Concerning the
(omitted) proofs, they are basically of two kinds: either
verifications via careful computations with certain special functions,
or geometric manipulations of (suitably decorated) triangulations of
$3$-manifolds. For the second ones, as well as for many properties of
the special functions, one has, in fact, to repeat with mild
modifications the arguments that are explained in full details in
$\{$Sections 1-6 $+$ the Appendix$\}$ of \cite{BB1}. The present paper
also accomplishes the `ideology' underlying Sections 7-9 of
\cite{BB1}, which are, by the way, completely surpassed. The results of
the present paper, limited to characters with values
in a Borel subgroup of $PSL(2,\mc)$, had been
announced in \cite{BB3}.

\medskip
 
\noindent In Section \ref{sumideal}, for any compact closed oriented 3-manifold
$W$ equipped with a  character $\rho$ with values in $PSL(2,\C)$, one
introduces special descriptions of $(W,\rho)$ via simplicial
$PSL(2,\C)$-valued 1-cocycles, called $\Dd$-{triangulations}.  A
simple and geometrically meaningful procedure of {\it idealization}
converts every $\Dd$-triangulation $\Tt$ into an $\Ii$-triangulation
$\Tt_{\Ii}$ made by hyperbolic ideal tetrahedra with ordered
vertices. 

\smallskip

\noindent In Section \ref{sumdilog}, one constructs a suitably uniformized version R of the classical {\it Rogers dilogarithm} L. We interpret R as a function of ideal tetrahedra with ordered vertices, enriched with suitable
$\Z$-valued decorations of the edges that induce {\it log-branches} on
the triples of moduli. It behaves well up to
tetrahedral symmetries, and it verifies nice functional five term
identities which lift arbitrary $2\leftrightarrow 3$ moves on ideal
triangulations. This section is largely inspired by W. Neumann's papers \cite{N1,N2}.

\smallskip

\noindent In Section \ref{sumdinv},  one states that the signed sum ${\rm R}(\Tt_{\Ii},f)$ of
the R-values of the tetrahedra of any fixed idealized triangulation
$\Tt_{\Ii}$ (suitably enriched with a global {\it flattening} $f$, which
induces a coherent log-branch system) well defines an invariant
${\rm R}(W,\rho)\in \C/(\pi^2/6)\Z$, called the {\it dilogarithmic
invariant}. This refines the $\C/\pi^2\Q$-valued Dupont-Sah's
description of the second Cheeger-Chern-Simons class for flat $SL(2,\mc)$-bundles. 
In particular, ${\rm R}(W,\rho)$ recovers the volume and the Chern-Simons invariant of the character
$\rho$. We stress that the proof of invariance does not use any 
sophisticated group-cohomological argument: one shows directly that any 
two arbitrary (enriched) $\Ii$-triangulations for $(W,\rho)$ can be connected by so-called (enriched) $\Ii$-{\it transits}.
 These transits are supported by usual elementary moves on bare triangulations, 
and dominated by $\Dd$-{\it transits} between $\Dd$-triangulations.

\smallskip

\noindent The dilogarithmic invariant emerged as a natural `classical'
counterpart from the effort in understanding the structure of the
{\it quantum hyperbolic invariants} (QHI), that we had 
constructed in $\{$Sections 1-6 $+$ the Appendix$\}$ of \cite{BB1} 
for triples $(W,L,\rho)$, where $L$ is a {\it non-empty}
link in $W$ and $\rho$ is only a $B$-valued character ($B$ is the
Borel subgroup of $SL(2,\C)$ of upper triangular matrices). The main
ingredient of the QHI are the Faddeev-Kashaev (non symmetric) {\it
quantum dilogarithms}. In Section \ref{sumQHI}, one revisits the $B$-QHI construction in the spirit of the dilogarithmic invariant's one, and one realizes that:

\smallskip

(1) For any $PSL(2,\C)$-character $\rho$, by using the quantum
dilogarithms one can define state sums supported by the idealization
$\Tt_{\Ii}$ of any $\Dd$-triangulation $\Tt$ for $(W,\rho)$. These
state sums do not yet define an invariant because they do not behave
well up to tetrahedral symmetries, and they verify only some special
instances of five term identities.

(2) A specific procedure of {\it local} symmetrization of the quantum
dilogarithms leads to fix an arbitrary non-empty link $L$ in $W$ in
order to fix one coherent globalization. This is supported by any
$\Dd$-triangulation $\Tt$ for $(W,\rho)$ in which the link $L$ is
realized as a Hamiltonian subcomplex $H$. The symmetrization of the
quantum dilogarithms is governed by any fixed {\it integral charge} $c$
relative to $H$. A charge is formally similar to the flattenings
used for the dilogarithmic invariant, the main difference being that
it is moduli-independent. It has a global structure which depends on
the link-fixing, and eventually encodes the link itself.  All this
gives the notion of $\Dd$-triangulation $(\Tt,H,c)$ for the triple
$(W,L,\rho)$, and any such triangulation supports suitably modified
state sums $H_N(\Tt_{\Ii},c)$ indexed by the odd integers $N>1$.

(3) The modified state sums satisfy all instances of five term identities,
so they eventually well define (up to a sign and a N-th root of unity
multiplicative factor) the QHI $H_N(W,L,\rho)$ for pairs $(W,L)$
equipped with arbitrary $PSL(2,\C)$-characters.
\smallskip

\smallskip

\noindent Every state sum $H_N(\Tt_{\Ii},c)$ looks formally very like $\exp
(R(\Tt_{\Ii},f))$, with the charge $c$ playing the role of the
flattening $f$. The presence  of the link in the construction of
$H_N(W,L,\rho)$ is entirely a consequence of the adopted specific
symmetrization procedure of the quantum dilogarithms. In order to get the general QHI for arbitrary $PSL(2,\C)$-characters, no further `quantum algebra' than for the $B$-QHI is necessary. In fact the
old $B$-state sums slightly differ from $H_N(\Tt_{\Ii},c)$ 
by a scalar factor which is expressed
in terms of $(\Tt,c)$, not only of $(\Tt_{\Ii},c)$. 

\smallskip

\noindent In Section \ref{sumscissors} one places the dilogarithmic invariant
and the QHI in the framework of the theory of {\it scissors congruence classes}.

\smallskip

\noindent One expects that the dilogarithmic invariant plays a main role to express the (dominant term of the) asymptotic expansion of the QHI when $N\to
\infty$. This {\it Volume Conjecture} for $PSL(2,\C)$-QHI is discussed 
in Section \ref{sumvol} where it is also supported by some results about the actual asymptotic behaviour of the quantum dilogarithms.

\smallskip

\noindent In Section \ref{sumcusped} one considers the extension of the
dilogarithmic invariant and of the QHI to {\it cusped} hyperbolic
3-manifolds. This leads us to discuss the relationship between the
corresponding Volume Conjecture and a current Volume Conjecture for
the colored Jones invariant $J_N(L)$ of hyperbolic knots $L$ in $S^3$.

\medskip

\noindent Finally we stress that our technic of `straightening' (generic)
1-cocycles, which applies to $(PSL(2,\C), \partial \bar{\mh}^3)$ in the
present paper and produces the idealization, looks very general. It
can be considered for other pairs $(G,X)$, where $G$
is a Lie group acting transitively on $X$ and there is a suitable
notion of {\it straight} tetrahedra with vertices in $X$. Cases of
most interest are $(SU(2), S^3)$, $(SO(3),S^3)$,
$(PSL(2,\C),\bar{\mh}^3)$. Moreover it should be interesting to extend it
to a treatment of  `truncated tetrahedra' (hence manifolds with
boundary).
By using other 'potential functions' one could
obtain the analogue of the dilogarithmic invariant for $W$ endowed
with $G$-valued characters. This could help to unterstand
(and possibly refine) other state sum quantum-invariants, such as the
Turaev-Viro one, in terms of scissors congruence classes (and, in particular,
give a geometric interpretation to
their asymptotic behaviour). We plan to study these generalizations in
future works.
\medskip

%% file: SUMIDEAL.tex
\section{$\Dd$-triangulations and the idealization} \label{sumideal}
Let $W$ be a compact, closed, oriented $3$-manifold, and $\rho$ be a
flat principal bundle over $W$ with structural group $PSL(2,\C)$.  We
consider $\rho$ up to isomorphisms of flat bundles; equivalently, it
can be identified with a conjugacy class of representations in
$PSL(2,\C)$ of the fundamental group of $W$, i.e. with a
$PSL(2,\C)$-character of $W$.  The pairs $(W,\rho)$ are considered up
to oriented homeomorphisms. By using the hauptvermutung, depending on
the context, we will freely assume that $W$ is endowed with a
(necessarily unique) PL or smooth structure, and use differentiable or
PL homeomorphisms.

\medskip

\noindent 
Given any singular triangulation $T$ of $W$ (i.e. a triangulation
where tetrahedra may have self and multiple adjacencies), and any
system of orientations of the edges of $T$, $\rho$ can be represented
by $PSL(2,\C)$-valued simplicial $1$-cocycles, up to coboundaries of
$0$-cochains. For our purposes, we need to specialize the kind of
triangulations, orientations and $1$-cocycles.

\noindent 
Any $T$ as above can be considered as a finite family $\{\Delta_i\}$
of {\it oriented abstract} tetrahedra, each being endowed with the
standard triangulation with $4$ vertices and the orientation induced
by the one of $W$, together with a system of identifications of pairs
of distinct (abstract) $2$-faces. We will often distinguish between
vertices, edges, 2-faces {\it in $T$}, that is after the
identifications, and {\it abstract} edges,.., that is of the abstract
$\Delta_i$'s. We view each $\Delta_i$ as positively embedded as a
straight tetrahedron in $\R^3$ endowed with the orientation specified
by the standard basis. Later we shall use the hyperbolic space
$\mh^3$. It is also oriented by stipulating that, in the disk model,
it is oriented as an open set of $\R^3$. The natural boundary
$\partial \bar{\mh}^3=\mc\mathbb{P}^1 = \C \cup\{ \infty \}$ of $\mh^3$ is
oriented by its complex structure. The action of $PSL(2,\C)$ on
$\mh^3$ and $\mc\mathbb{P}^1$ is the one of orientation preserving
isometries of $\mh^3$.
 
\begin{defi}\label{Dtriang} 
{\rm A $\Dd$-{\it triangulation} for the pair $(W,\rho)$ consists of a
triple $\Tt=(T,b,z)$ where:

\smallskip

(i) $T$ is a {\it quasi-regular} triangulation of $W$, that is every
edge of $T$ has two distinct vertices as endpoints;

\smallskip

(ii) $b$ is a {\it branching} of $T$, that is a system of orientations
of the edges such that the one induced on each abstract $\Delta_i$ is
associated to a total ordering $v_0,v_1,v_2,v_3$ of its (abstract) vertices: 
each edge is oriented by the arrow emanating from the
smallest endpoint;

\smallskip

(iii) $z$ is a $1$-cocycle on $(T,b)$ representing $\rho$ such that $(T,b,z)$ is {\it idealizable} 
(see Def. \ref{idealizable}).} 
\end{defi}

\begin{remark}\label{totalord} 
{\rm Quasi-regular (even regular indeed) triangulations of $W$ do
exist.  A total ordering of the vertices of a quasi-regular
triangulation $T$ clearly induces a branching. {\it In the present
paper we will only consider these special branchings associated to
total orderings of the vertices}. This simplifies certain proofs, but
all the results eventually hold true also for arbitrary
branchings. For more information about branchings, see \cite{BB1} and
\cite{BP2}}.
\end{remark}

\noindent 
Given a branching $b$ on a oriented tetrahedron $\Delta$ (realized in
$\R^3$ as above), one can define an orientation of any of its
simplices, not only of the edges. Denote by $E(\Delta)$ the set of
$b$-oriented edges of $\Delta$, and by $e'$ the edge opposite to
$e$. We put $e_0=[v_0,v_1]$, $e_1=[v_1,v_2]$ and
$e_2=[v_0,v_2]=-[v_2,v_0]$. The ordered triple of edges
$$(e_0=[v_0,v_1],e_2=[v_0,v_2],e_1'=[v_0,v_3])$$ departing from $v_0$
defines a \emph{$b$-orientation} of $\Delta$. We say that $(\Delta,b)$
is positive if its $b$-orientation agrees with the one of $\R^3$, and
negative otherwise; we indicate it by a sign $*=*_b=\pm 1$.  The
$2$-faces of $\Delta$ can be named by their opposite vertices. We
orient them by working as above on the boundary of each 2-face $f$:
there is a $b$-ordering of the vertices of $f$, and an orientation of
$f$ which induces on $\partial f$ the prevailing orientation among the
three $b$-oriented edges. If $z$ is $PSL(2,\C)$-valued 1-cocycle on
$(\Delta,b)$, we write $z_j = z(e_j)$ and $z'_j = z(e'_j)$. Then, one reads
the cocycle condition on the $2$-face opposite to $v_3$ as $ z_0z_1z_2^{-1} = 1$. 
These considerations
apply to each abstract tetrahedron of any branched triangulation
$(T,b)$ of $W$ and to (the restrictions of) any $PSL(2,\C)$-valued
1-cocycle $z$ on $(T,b)$.

\begin{defi}\label{idealizable}
{\rm Let $(\Delta,b,z)$ be a branched tetrahedron endowed with a
$PSL(2,\C)$-valued 1-cocycle $z$. It is {\it idealizable} iff
$$ u_0=0,\ u_1= z_0(0),\ u_2= z_0z_1(0),\ u_3= z_0z_1z'_0(0) $$ are 4 
distinct points in $\C \subset \mc\mathbb{P}^1= \partial \bar{\mh}^3$
which span an non degenerate hyperbolic ideal tetrahedron (with ordered
vertices). A triangulation $(T,b,z)$ is {\it idealizable} iff all its
tetrahedra  $(\Delta_i,b_i,z_i)$ are idealizable.}
\end{defi} 

\noindent 
If $(\Delta,b,z)$ is idealizable, for all $j=0,1,2$ one can associate
to $e_j$ and $e'_j$ the same {\it cross-ratio} modulus $w_j\in
\C\setminus \{0,1\}$ of the hyperbolic ideal tetrahedron defined by $(u_0,u_1,u_2,u_3)$. 
By using the cyclic ordering of the edges induced by $b$, one has
$$w_{j+1} = 1/(1-w_j)$$
and
$$ w_0 = ( u_2- u_1)u_3/u_2(u_3- u_1) := -p_1/p_2\ .$$ 
Similarly we shall write 
$$w_j = - p_{j+1}/p_{j+2}\ .$$ 
\noindent Set $w=(w_0,w_1,w_2)$ and call it a {\it modular triple}. As the ideal tetrahedron is non-degenerate, 
the imaginary parts of the $w_j$'s  are not equal to zero, and share the same sign $*_w = \pm 1$.

\begin{defi}\label{idealizations}{\rm  We call $(\Delta,b,w)$ the {\it idealization} of the idealizable $(\Delta,b,z)$.  
For any
 $\Dd$-triangulation $\Tt=(T,b,z)$ of $(W,\rho)$, its {\it
 idealization} $\Tt_{\Ii}=(T,b,w)$ is given by the family $\{(\Delta_i, b_i, w_i)\}$ of idealizations of the 
$(\Delta_i,b_i, z_i)$'s. We say that $\Tt_{\Ii}$ is an $\Ii$-triangulation.}
\end{defi}

\begin{remarks}\label{idealrem}{\rm (1) We have inglobed the `non-degenerate'
assumption into the notion of `idealizable'. This simplifies the exposition and 
also certain proofs; however, it is not necessary for the validity of the results of the paper.

\smallskip

\noindent (2) If $z_j$ acts on $\mc\mathbb{P}^1$ as $(a_jx +b_j)/(c_jx+d_j)$, then $z_j(0) = b_j/d_j$. 
Thus it is immediate to formulate the `idealizability' condition in terms of a simple system of real algebraic 
inequalities on the entries of the $z_j$'s.

\smallskip

\noindent (3) In \cite{BB1} we have used so-called \emph{full} $B$-valued 1-cocycles $z$ to construct the $B$-QHI.  
This means that for any
edge $e$ the upper-diagonal entry $x(e)$ of $z(e)$ is non-zero. It is
easy to verify that a $B$-cocycle is full iff it is idealizable
(forgetting the non-degenerate assumption). In
\cite{BB3} we proposed an idealization of full $B$-cocycles, which was indeed
a specialization of the present general procedure. 
The idealization of a $\Dd$-tetrahedron with a full $B$-cocycle is simply given 
by $ w_j = - q_{j+1}/q_{j+2}$, where $q_j=x(e_j)x(e_j')$ for $j=0,1$, and $q_2=-x(e_2)x(e_2')$.

\smallskip

\noindent (4) It follows from the cocycle condition that $p_0 + p_1 + p_2 =0$
(also that $q_0 + q_1 + q_2 =0$ for full $B$-cocycles).}
\end{remarks}

\noindent The following lemma is immediate:

\begin{lem} \label{esitenza} For any $PSL(2,\C)$-character $\rho$, any 
quasi-regular branched triangulation $(T,b)$ of $W$ can be completed to a $\Dd$-triangulation $(T,b,z)$ 
of the pair $(W,\rho)$.
\end{lem}

\noindent In fact, given any 1-cocycle, one can perturb it by the coboundaries of generic $0$-cochains which 
are injective on the vertices of $T$, and get an idealizable one.

\medskip

\noindent{\bf Tetrahedral symmetries.} The following simple lemma states the good behaviour of 
the idealization with respect to a change of branching (i.e. the `tetrahedral symmetries').

\begin{lem} \label{Iisym}  Denote by $S_4$ the permutation group on four elements.
A permutation $p \in S_4$ of the vertices of an idealizable tetrahedron 
$(\Delta,b,z)$ gives another idealizable tetrahedron $(\Delta,b',z')$. 
For every edge $e$, $z'(e)=z(e)^{-1}$ iff the orientations of $e$ for $b$ and $b'$ are opposite; 
otherwise  $z'(e)=z(e)$. The permutation turns the idealization $(\Delta,b,w)$ into $(\Delta,b',w')$, where, 
for each edge $e$, 
$w'(e)=w(e)^{\epsilon(p)}$, and $\epsilon(p)$ is the signature of $p$.  
\end{lem}

\noindent Consider for instance the transposition $(0,1)$. It turns the set of 4 points
$$ 0, z_0(0), z_0z_1(0), z_0z_1z'_0(0) $$
into
$$ 0, (z_0)^{-1}(0), z_1(0), z_1z'_0(0) \ .$$
By applying on the second set of points the hyperbolic isometry $z_0$ one gets the initial 
set after the transposition of its first two members, and similarly for other permutations. 
The lemma follows immediately, due to the behaviour of the cross-ratio up to vertex permutation.

\begin{remark}{\rm As the signature of a permutation also changes the sign $*_b$ of the branching, 
the above lemma is coherent with the usual symmetry relations holding in the classical (pre)-Bloch group:
$[x] = [1/(1-x)]$ and $[x]=-[1/x]$ (see e.g. \cite{D2}, \cite{N1}).}
\end{remark}

\smallskip

\noindent {\bf Hyperbolic edge compatibility.} We are now concerned with an important global property of the 
idealized triangulations $\Tt_{\Ii}$. Before to state it, let us stress that when dealing with modular triples one has 
to be careful with the orientations. Recall that every $\Ii$-tetrahedron
$(\Delta,b,w)$ is oriented by definition; in the case of an $\Ii$-triangulation this is given by the orientation of $W$. 
There is also the $b$-orientation encoded by the sign $*=*_b$. The idealization `physically' realizes the vertices 
of $\Delta$ on $\partial \bar{\mh}^3$, with the ordering induced by $b$. So the $b$-orientation may or may not agree 
with the orientation of the spanned ideal tetrahedron induced by the fixed orientation of $\mh^3$, which is 
encoded by the sign $*_w$ of the modular triple.

\begin{defi}\label{geometrico}{\rm  We say that $(\Delta,b,w)$ is {\it geometric} iff $*_b = *_w$.}
\end{defi}

\noindent Given any $\Ii$-triangulation $\Tt_{\Ii}=(T,b,w)$, the contribution of each $(\Delta_i,b_i,w_i)$ to 
any computation with the moduli is given by the $w(e)^{*}$'s, where $e$ is any edge in $\Delta_i$ and $*= *_{b_i}$. 
The next Lemma \ref{compat} is a first concretization of this fact (see also the notion of $\Ii$-transit below).
\medskip

\noindent  Let $\Tt_{\Ii}=(T,b,w)$ be as above. 
Denote by $E(T)$ the set of edges of $T$,  by $E_{\Delta}(T)$ the whole set of edges
of the associated abstract tetrahedra $\{\Delta_i\}$, and by $\epsilon:E_{\Delta}(T) \longrightarrow E(T)$ 
the natural identification map. 


\begin{lem} \label{compat} For any edge $e \in E(T)$, we have 
$\textstyle \prod_{a \in \epsilon^{-1}(e)} w(a)^{*} =1$,  
where $*=\pm 1$ according to the $b$-orientation of the tetrahedron $\Delta_i$ that contains $a$.
\end{lem}

\noindent 
This means that the $w(a)^{*}$'s around each $e$ verify the usual
compatibility condition, necessary when one
tries to construct hyperbolic $3$-manifolds by glueing ideal
tetrahedra. A key point in the proof is that around every edge $e$
of $T$ one may only meet an even number of tetrahedra such that the
$b$-orientations of the two faces containing $e$ are opposite. This
follows from the fact that $W$ is orientable.

\begin{remark}\label{geometricD}{\rm  
Every idealized triangulation $\Tt_{\Ii}$ for $(W,\rho)$ necessarily
includes some {\it non}-geometric ideal tetrahedra.  In fact if all
the tetrahedra were geometric, the compatibility condition in Lemma
\ref{compat} would imply that the sum of the arguments is exactly
equal to $2\pi$ around all the edges of $T$, not only an even multiple
of $2\pi$ (see Lemma E.6.1. of \cite{BP1}). So, one should have a
genuine hyperbolic manifold structure on $W$ out of the vertices of
$T$, with $W$ being triangulated by embedded geodesic ideal
tetrahedra. Using tetrahedra truncated by suitable horospheres, one
deduces that the spherical link surrounding each vertex of $T$ would
inherit a $(\C^*,\mc)$-structure, which is impossible.}
\end{remark}
\smallskip

\begin{figure}[ht]
\begin{center}
\includegraphics[width=7cm]{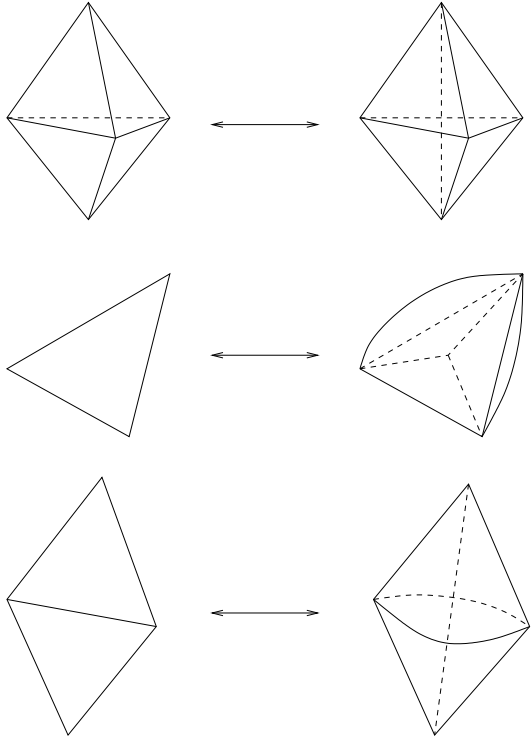}
\caption{\label{figmove1} the moves between singular triangulations.} 
\end{center}
\end{figure}

\noindent {\bf $\Dd$- and $\Ii$-transits.} 
It is well-known that given any two (singular) triangulations $T_0$
and $T_1$ of $W$ there exists a finite sequence $T_0 \rightarrow
\ldots \rightarrow T_1$ of so-called $2 \leftrightarrow 3$ and {\it
bubble} {\it moves} that turns $T_0$ into $T_1$. For instance, it is a
consequence of the duality between ideal triangulations and standard
spines of $T \setminus T^0$ (removing the vertices), and the calculus
for standard spines due to Matveev \cite{Mat} and Piergallini
\cite{Pi}. It is sometimes convenient, for technical reasons (for
instance when dealing with arbitrary branchings as in \cite{BB1}), to
consider a further so-called $0 \leftrightarrow 2$ {\it move}. See
Fig. \ref{figmove1}.

\noindent Next we consider moves on $\Dd$-triangulations $\Tt=(T,b,z)$ and $\Ii$-triangulations 
$\Tt_{\Ii}=(T,b,w) $ for the pair $(W,\rho)$, called $\Dd$-{\it transits} and $\Ii$-{\it transits} respectively. 
They are supported by the above bare triangulation
moves, but they also include the transits of the respective extra-structures.

\noindent First of all one requires that the condition to be quasi-regular is preserved by the moves.
We stress that this is not an automatic fact, on the contrary this leads
to one main technical complication in the proofs.

\noindent A move $T_0 \leftrightarrow T_1$ between branched triangulations $(T_0,b_0)$ and $(T_1,b_1)$ 
defines a {\it branching transit} $(T_0,b_0) \leftrightarrow
(T_1,b_1)$ if $b_0$ and $b_1$ agree on the common edges.  As
stipulated above, we shall only consider branched quasi-regular
triangulations, where the branchings are defined by total orderings of
the vertices. In such a case any $2 \leftrightarrow 3$ or $0
\leftrightarrow 2$ move that preserves the quasi-regularity of the
triangulations can be completed in a unique way to a branching
transit, but a bubble move may be completed in different ways, each of
them being a possible transit.
 
\noindent Let $(T_0,b_0)$, $(T_1,b_1)$ be as above and $z_k \in Z^1(T_k;PSL(2,\C))$, $k=0,1$. We have a 
(resp. {\it idealizable}) {\it cocycle transit} $(T_0,b_0,z_0)
\leftrightarrow (T_1,b_1,z_1)$ if $z_0$ and $z_1$ agree on the common
edges (resp. and both are idealizable $1$-cocycles). Note that for $2
\rightarrow 3$ and $0 \rightarrow 2$ moves, given $z_k$ there is only
one (resp. at most one) $z_{k+1}$ with this property. We stress that
in some special cases a $2\to 3$ transit of an idealizable cocycle can
actually not preserve the idealizability, but generically this does
not hold.  For positive bubble moves there is always an infinite set
of possible (idealizable) cocycle transits.

\noindent We say that $(T_0,b_0,z_0) \leftrightarrow (T_1,b_1,z_1)$ as above is a 
$\Dd$-transit when both $z_0$ and $z_1$ are idealizable.

\medskip

\noindent Let us now consider the transit for the idealized triangulations.
Consider the convex hull of five distinct points 
$u_0,u_1,u_2,u_3,u_4 \in \partial \bar{\mathbb{H}}^3$, with the two possible triangulations $Q_0$ $Q_1$ 
made of the oriented hyperbolic ideal tetrahedra $\Delta^i$ obtained by omitting $u_i$. An edge $e$ of $Q_i 
\cap Q_{i+1}$ belongs to one tetrahedron of $Q_i$ iff it belongs to two tetrahedra of $Q_{i+1}$. 
Then, the modulus 
of $e$ in $Q_i$ is the product of the two moduli of $e$ in
$Q_{i+1}$. Also, the product of the moduli on the central edge of
$Q_1$ is equal to $1$.

\noindent Let $T \to T'$ be a $2\to 3$ move. 
Consider the two (resp. three) abstract 
tetrahedra of $T$ (resp. $T'$) involved in the move. They determine subsets 
$\widetilde{E}(T)$ of $E_{\Delta}(T)$ 
and $\widetilde{E} (T')$ of $E_{\Delta}(T')$. Denote their complementary sets by $\widehat{E} (T)$ 
and $\widehat{E} (T')$. Clearly one can identify $\widehat{E} (T)$ and $\widehat{E} (T')$. Using the above 
configurations $Q_0$ and $Q_1$, and recalling the considerations 
made before Lemma \ref{compat}, we are led to the following definition:

\begin{figure}[ht]
\begin{center}
\includegraphics[width=11cm]{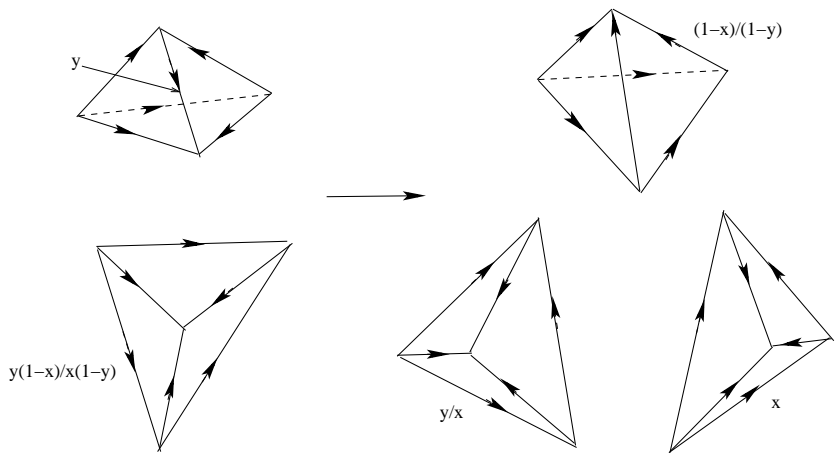}
\caption{\label{idealt} a $2 \leftrightarrow 3$ ideal transit.} 
\end{center}
\end{figure}

\begin{defi} \label{id-transit} 
{\rm One has a $2 \to 3 \ $ {\it $\Ii$-transit} $(T,b,w) \rightarrow (T',b',w')$ of $\Ii$-triangulations for a pair $(W,\rho)$ 
if:  

\smallskip

\noindent 1) $w$ and $w'$ agree on $\widehat{E} (T) = \widehat{E} (T')$;

\noindent 2) for each common edge $e\in \epsilon_T(\widetilde{E} (T))\cap \epsilon_{T'}(\widetilde{E} (T'))$ one has}
\begin{equation}\label{ideqmod} \prod_{a\in \epsilon_T^{-1}(e)}w (a)^*=
 \prod_{a'\in \epsilon_{T'}^{-1}(e)}w' (a')^*\ ,
\end{equation}
\noindent {\rm where $*=\pm 1$ according to the $b$-orientation of the tetrahedron that contains $a$ (resp. $a'$). 
One has a $0 \to 2$ (resp. \emph{bubble}) \emph{$\Ii$-transit} if the above first condition is satisfied, and 
one replaces the second by: 

\smallskip

\noindent 2') for each edge $e\in \epsilon_{T'}(\widetilde{E} (T'))$ one has}
\begin{equation}\label{ideq2mod} \prod_{a'\in \epsilon_{T'}^{-1}(e)}w' (a')^*=1\ .
\end{equation}
\end{defi}  

\noindent $\Ii$-transits for negative $3 \to 2$ moves are defined in exactly the same way, 
and for negative $2 \to 0$ 
and bubble moves $w'$ is defined by simply forgetting the moduli of
the two disappearing tetrahedra. The condition (1) above, implies that the product of the $w'(a')^*$'s
around the new edge is equal to $1$.  
A $2 \leftrightarrow 3$
$\Ii$-transit is shown in Fig. \ref{idealt}; we only indicate the
first component of each modular triple. In general, the relations
(\ref{ideqmod}) may imply that $w$ or $w'$ equals $0$ or $1$ on some
edges.  In that case, the $2 \leftrightarrow 3$ $\Ii$-transit
fails. In particular, in Fig. \ref{idealt} we assume that $x \ne y$.

\noindent Note that for $2 \leftrightarrow 3$ $\Ii$-transits $w'$ is uniquely determined by $w$, whereas for $0 \to 2$ 
and bubble $\Ii$-transits there is one degree of freedom in choosing $w'$. Also, condition (\ref{ideq2mod}) simply 
means that such transits give \emph{the same} modular triples to the two new tetrahedra, for their $b$-orientations
are opposite.

\begin{remark}\label{blochrel}{\rm Once it is expressed in terms of the
involved moduli $w_0$'s, any $2 \leftrightarrow 3$ $\Ii$-transit
dominates an instance of the {\it five term identities} which enter in
the construction of the classical (pre)-Bloch group. For example, the
transit in Fig. \ref{idealt} corresponds to the relation
$$[x]-[y]+[y/x]-[(1-x^{-1})/(1-y^{-1})]+[(1-x)/(1-y)]=0\ .$$}
\end{remark}
The next proposition states the remarkable fact that $\Dd$-transits
and $\Ii$-transits together with the idealization make commutative
diagrams, that is the $\Dd$-transits {\it dominate} the
$\Ii$-transits.
\begin{prop}\label{DdomI}  
Let $\mathfrak{d}$ be any $\Dd$-transit of $\Dd$-triangulations 
and $\mathfrak{i}$ be any $\Ii$-transit of $\Ii$-triangulations
for a pair $(W,\rho)$. Denote by $\Ii$
the idealization map $\Tt \to \Tt_{\Ii}$. Then for every
$\mathfrak{d}$ there exists $\mathfrak{i}$ (resp. for every
$\mathfrak{i}$ there exists $\mathfrak{d}$) such that $\mathfrak{i}
\circ \Ii = \Ii \circ \mathfrak{d}$.
\end{prop}
For $2 \leftrightarrow 3$ transits there is also an uniqueness
statement.  The proof is not hard. By using the tetrahedral symmetries
of Lemma \ref{Iisym}, it is enough to show the proposition for one branching
transit configuration (for instance the one of Fig.
\ref{idealt}). The verification follows almost immediately from the
definition of the idealization, as for Lemma \ref{Iisym}.  Note that
the possible failures of $2\to 3$ transits of idealizable cocycles
that we mentionned above correspond to the failures of $2\to 3$
$\Ii$-transits (for instance when $x=y$ in Fig.  \ref{idealt}).

%% file: NSUMDILOG.tex
\section{Classical dilogarithms}\label{sumdilog}

\subsection{Rogers dilogarithm}\label{rog}

\noindent
Denote by $\log$ the standard branch of the logarithm, with arguments
in $]-\pi,\pi]$.  Put $\mathfrak{D} = \mc \setminus \{(-\infty;0) \cup
(1;+\infty) \}$.  The \emph{Rogers dilogarithm} is the complex
analytic function defined over $\mathfrak{D}$ by
\begin{equation}\label{Rdilog}
 {\rm L}(x) = -\frac{\pi^2}{6} -\frac{1}{2} \int_0^x \biggl( \frac{\log(t)}{1-t} + \frac{\log(1-t)}{t} \biggr) \ dt \ ,
\end{equation}
\noindent where we integrate first along the path $[0;1/2]$ on the real axis and then 
along any path in $\mathfrak{D}$ from $1/2$ to $x$. Here we add $-\pi^2/6$ so that ${\rm L}(1)=0$. 
For $\vert x \vert <1$, one may also write L as 
$$ {\rm L}(x) = -\frac{\pi^2}{6} + \frac{1}{2} \log(x)\log(1-x) + \sum_{n=1}^{\infty} \frac{x^n}{n^2}\ .$$
 \noindent The sum in the right-hand side is the power series expansion in the open unit disk of the 
\emph{Euler dilogarithm} ${\rm Li}_2$, defined over $\mc \setminus (1;+\infty)$ by 
$${\rm Li}_2(x) =  - \int_0^x \frac{\log(t)}{1-t} \ dt\ .$$
\noindent For a detailed study of the dilogarithm functions and their relatives, see \cite{L} or the review \cite{Z}. Below we will also call `Rogers dilogarithm' the multi-valued function obtained from L via analytic continuation. The function L is related to the \emph{Bloch-Wigner dilogarithm}
$${\rm D}_2(x) = {\rm Im} \bigl( {\rm Li}_2(x) \bigr) + \arg(1-x)\log\vert x\vert\ ,$$
\noindent which is obtained by adding to ${\rm Im}( {\rm Li}_2(x))$ the needed correction term to compensate its 
jump along the branch cut $(1;+\infty)$. The function ${\rm D}_2(x)$ is a real analytic continuation of 
${\rm Im}( {\rm Li}_2(x))$ on $\mc \setminus \{0,1\}$, and it is continuous (but not differentiable) at $0$ and $1$.
 It gives the volume of oriented hyperbolic tetrahedra by the formula
\begin{equation}\label{BWvol}
{\rm Vol}(\Delta,b,w) = {\rm D}_2(w_0)\ ,
\end{equation}
\noindent where, with the notations of Section \ref{sumideal},
$w$ denotes the 
modular triple with respect to a positive branching $b$.
In particular, we have the 6-fold symmetries
\begin{equation}\label{BWsym}
{\rm D}_2(w_0) = {\rm D}_2(w_1)= {\rm D}_2(w_2)=- {\rm D}_2(w_0^{-1}) = -{\rm D}_2(w_1^{-1}) = 
-{\rm D}_2(w_2^{-1})\
.\end{equation}
\noindent Moreover, if we apply the formula (\ref{BWvol}) to Fig. \ref{idealt} we get the five term functional relation
\begin{equation}\label{BWfivet}
 {\rm D}_2(y) + {\rm D}_2(\frac{1-x^{-1}}{1-y^{-1}}) = {\rm D}_2(x)
 +{\rm D}_2(y/x) + {\rm D}_2(\frac{1-x}{1-y})
\end{equation} 
\noindent when $x \ne y$. Finally, all the other five term relations obtained by changing the branching in 
Fig. \ref{idealt} also hold true, due to (\ref{BWsym}). 

\noindent One would like to think of the Rogers 
dilogarithm L as the natural complex analytic analogue of ${\rm D}_2(x)$. But L verifies similar 
five term relations only by putting restrictions on the variables. For instance, the analog of 
(\ref{BWfivet}) is the so-called Schaeffer's identity 
\begin{equation}\label{Rfivet}
{\rm L}(x) - {\rm L}(y) +{\rm L}(y/x) - {\rm L}(\frac{1-x^{-1}}{1-y^{-1}}) + {\rm L}(\frac{1-x}{1-y})=0
\end{equation}
which for real $x$, $y$ holds only when $0 < y < x < 1$. In fact, this
identity characterizes the Rogers dilogarithm: if $f(0;1) \to \mr$ is
a 3 times differentiable function satisfying (\ref{Rfivet}) for all $0
<y <x <1$, then $f(x) = k{\rm L}(x)$ for a suitable constant $k$
\cite[Sect. 4]{Ro}, \cite[App.]{D1}.  By analytic continuation, the
relation (\ref{Rfivet}) holds true for complex parameters $x$, $y$,
providing that the imaginary part of $y$ is $\geq 0$, and $x$ lies
inside the triangle formed by $0$, $1$ and $y$. For such $x$, $y$,
note that also all the other arguments of L in (\ref{Rfivet}) have
positive imaginary parts. With this restriction, this relation
corresponds to one specific instance of $\Ii$-transit (see Remark
\ref{blochrel}). This naively suggests the possibility to set 
$${\rm L}(\Delta,b,w) = {\rm L}(w_0)\quad ,$$ 
and try to use it to build an invariant for $(W,\rho)$ that should be computable by using any idealized
triangulation $\Tt_{\Ii}$. However, one is
immediately faced to the following correlated difficulties:

\medskip

(1) ({\it Uniformization}) One has to deal with the different branches of L.
\smallskip

(2) ({\it Symmetrization }) One realizes that $ {\rm L}(\Delta,b,w)$  only respects the
tetrahedral symmetries up to some elementary functions (see below).  
\smallskip

(3) ({\it Completing the five term relations}) As remarked above, L satisfies certain five term relations which correspond only to some peculiar instances of  $\Ii$-transits, and with restrictions on the range of moduli.
\medskip

\noindent In the next three subsections we outline how to solve these difficulties.
\subsection{Uniformization}\label{unif}
\noindent Let $\widehat{\mc} = \widehat{\mc}_{00} \cup \widehat{\mc}_{01} \cup \widehat{\mc}_{10} 
\cup \widehat{\mc}_{11}$, where $\widehat{\mc}_{\varepsilon \varepsilon'}$ ($\varepsilon, \varepsilon'=0,1$) is 
the Riemann surface of the function defined on $\mathfrak{D}=\mc
\setminus \{(-\infty;0) \cup (1;+\infty) \}$ by
$$x \mapsto (\log(x) +\varepsilon i\pi ,\log((1-x)^{-1})+\varepsilon'i\pi )\ .$$
\noindent Thus $\widehat{\mc}$ is the abelian cover of $\mc \setminus \{0,1\}$ obtained from $\mathfrak{D} \times 
\mz^2$ by the identifications
$$\begin{array}{l}
\{(-\infty;0)+i0\} \times \{p\} \times \{q\} \sim \{(-\infty;0) -i0\} \times \{p+2\} \times \{q\}\\
\{(1;+\infty)+i0\} \times \{p\} \times \{q\} \sim \{(1;+\infty)-i0\} \times \{p\} \times \{q+2\} \ ,
\end{array}$$
\noindent and the function 
\begin{equation}\label{logb1}
{\rm l}(x;p,q) = (\log(x) + pi\pi,\log((1-x)^{-1}) +qi\pi)
\end{equation}
\noindent is well-defined and analytic on $\widehat{\mc}$. Here $(-\infty;0)\pm i0$ comes from the upper/lower 
fold of $\mathfrak{D}$ with respect to $(-\infty;0)$, and similarly for $(1;+\infty)\pm i0$. Following \cite{N1,N2}, 
consider the following lift of the Rogers dilogarithm on $\widehat{\mc}$:
\begin{equation}\label{Ndef}
{\rm R}(x;p,q) = {\rm L}(x) + \frac{i\pi}{2} (p\log(1-x) +q\log(x))\ .
\end{equation}

\begin{lem}\label{R} The above formula well defines an analytic map ${\rm R}: \widehat{\mc}\to \C/\pi^2\Z$.
\end{lem}

\noindent One can view R as a uniformization mod($\pi^2$) of L. We want to interpret R as a function of our $\Ii$-tetrahedra $(\Delta,b,w)$. In order to do that, it is natural to enrich the decoration by a $\Z$-valued function $f$ on the edges of $\Delta$ such that, for every edge, $f(e)=f(e')$.  As for
$w=(w_0,w_1,w_2)$, one can write $f=(f_0,f_1,f_2)$ with respect to
$b$. Then we set
$$ {\rm R}(\Delta,b,w,f) = {\rm R}(w_0;f_0,f_1)\ .$$ Next we indicate
under which condition on $f$ the function ${\rm R}(\Delta,b,w,f)$ respects the
tetrahedral symmetries.

\subsection{Tetrahedral symmetries}
Here is a crucial definition.
\begin{defi} \label{flattening}
{\rm Let $(\Delta,b,w,f)$ and $f=(f_0,f_1,f_2)$ be as above. Set} 
 $${\rm l}_j(b,w,f)=\log(w_j) + f_ji\pi \ ,$$
\noindent {\rm for $j=1$, $2$, $3$. We say that $(f_0,f_1,f_2)$ is a \emph{flattening} of $(\Delta,b,w)$ if} 
$${\rm l}_0(b,w,f) + {\rm l}_1(b,w,f) +{\rm l}_2(b,w,f) = 0\ .$$ {\rm In that case, we call ${\rm l}_j(b,w,f)$ a \emph{log-branch} 
of $(\Delta,b,w)$ for the edge $e_j$, and set ${\rm l}(b,w,f)=({\rm l}_0(b,w,f),{\rm l}_1(b,w,f),{\rm l}_2(b,w,f))$ for the total log-branch associated to $f$.}
\end{defi}

\noindent Note that if $f$ is a flattening of $(\Delta,b,w)$, then it is a flattening of $(\Delta,b,u)$ for
every modular triple $u$ sufficiently close to $w$.

\noindent Let $(\Delta,b',w,f)$ be any enriched $\Ii$-tetrahedron.
By acting with a permutation $p \in S_4$ on the vertices of $\Delta$,
one passes from $b'$ to a new branching $b$. Denote by
$(\Delta,b,w,f)$ the new tetrahedron, where one still associates to
every edge the same values $w(e)$ and $f(e)$ as before, but they are
renamed according to the new ordering of the vertices given by
$b$. Let $\epsilon(p)$ be the signature of $p$. Set $w^{\epsilon(p)}=
(w_0^{\epsilon(p)}, w_1^{\epsilon(p)},w_2^{\epsilon(p)})$ and
$\epsilon(p)f=(\epsilon(p)f_0,\epsilon(p)f_1,\epsilon(p)f_2)$. We have

\begin{lem}\label{Rsym3} For any enriched $\Ii$-tetrahedron $(\Delta,b,w,f)$ the identities  
$${\rm R}(\Delta,b',u,f) =
\epsilon(p)\ {\rm R}(\Delta,b,u^{\epsilon(p)},\epsilon(p)f))\quad {\rm mod}(\pi^2/6) \mz$$
hold true for every permutation $p$ and for every modular triple $u$
sufficiently close to $w$ if and only if $f$ is a flattening of $(\Delta,b,w)$.
\end{lem}
\subsection{Complete five term relations}

\noindent We first define a notion of transit between flattened tetrahedra, and then use it to get 
all the required five term relations, without restrictions neither on the
underlying branching transit configuration, nor on the range of the
moduli.

\noindent Consider a $2 \to 3$ $\Ii$-transit  $(T,b,w) \rightarrow (T',b',w')$ as in Fig. \ref{idealt}. 
Give a flattening to each 
tetrahedron of the initial configuration, and denote by l the corresponding log-branch function on $T$. 
Recall the definition of the map $\epsilon_T$ before Lemma \ref{compat}. The simple idea is just to
formally take the $\log$ of the $\Ii$-transit.

\begin{defi}\label{lbtransit} {\rm A map $f': E_{\Delta}(T') \longrightarrow \mz$ defines a $2 \to 3$ 
{\it log-branch transit} $(T,b,w,f) \rightarrow (T',b',w',f')$  if for each common edge $e \in T \cap T'$ one has the following relation between log-branches:}
\begin{equation}\label{ideq} \sum_{a\in \epsilon_T^{-1}(e)}*\ {\rm l}(a)=
 \sum_{a'\in \epsilon_{T'}^{-1}(e)}*\ {\rm l}'(a')\ ,
\end{equation}
\noindent {\rm where $*=\pm 1$ according to the $b$-orientation of the tetrahedron that contains $a$ (resp. $a'$). 
A map $f': E_{\Delta}(T') \longrightarrow \mz$ defines a $0 \to 2$ (resp. bubble) \emph{log-branch transit} if for 
each edge $e \in T'$ one has}
\begin{equation}\label{ideq2} \sum_{a'\in \epsilon_{T'}^{-1}(e)} *\ {\rm l}'(a')=0\ .
\end{equation}
\end{defi}  
\noindent One easily verifies that log-branch transits actually define flattened tetrahedra, and that for a 
$2 \to 3$ log-branch transit the sum of values of l' about the new edge is always equal to zero. 
So log-branch transits for negative $3 \to 2$ moves are defined in exactly the same way, except that we also 
require that this last condition holds. For negative $2 \to 0$ and bubble moves the log-branch transits are defined 
by simply forgetting the log-branches of the two disappearing tetrahedra. 
The flattenings of a log-branch transit, associated to a given $\Ii$-transit $(T,b,w) \rightarrow (T',b',w')$, 
actually define a log-branch transit for every $\Ii$-transit $(T,b,u) \rightarrow (T',b',u')$ if $u$ (resp. $u'$)
is a modular triple sufficiently close to $w$ (resp. $w'$). 

\noindent Note that the relations (\ref{ideq2}) mean that the two new tetrahedra have \emph{the same} log-branches, 
for their $b$-orientations are always opposite. 

\begin{prop}\label{Rtransit} Let $(T,b,w,f) \rightarrow (T',b',w',f')$ be a log-branch transit. Then we have
\begin{equation}\label{Rtransit2}
\sum_{\Delta \subset T} *\ {\rm R}(\Delta,b,w,f) = \sum_{\Delta' \subset T'} *\ {\rm R}(\Delta',b',w',f')\quad
 {\rm mod}(\pi^2/6) \mz \quad ,
\end{equation} 
\noindent
where $*=\pm 1$ according to the $b$-orientation of $\Delta$
 (resp. $\Delta'$).  \end{prop} Thanks to the tetrahedral symmetries,
 it is enough to prove the proposition for one branching transit
 configuration.  For instance, for $2 \leftrightarrow 3$ transits one
 uses the one in Fig. \ref{idealt}, which underlies Schaeffer's
 identity (\ref{Rfivet}). The key and delicate point consists in
 realizing that the log-branch transit condition determines a proper
 {\it analytic} subset $\widehat{\mathfrak{G}}$ of
 $\widehat{\mc}^5$, on which the five term relation (\ref{Rtransit2})
 corresponds to an {\it analytic} relation. Then one verifies that this
 relation holds true on a non empty open subset of $\widehat{\mathfrak{G}}$,
as a consequence of Schaeffer's identity (\ref{Rfivet}) over its complex domain of validity. \
Then the result follows by using the analytic continuation principle.

\begin{remark}\label{duetoN}
{\rm The function ${\rm R}(x;p,q)$ with its five term relation that
lifts Schaeffer's identity was already considered by W. Neumann
\cite{N1,N2}. A peculiarity of our treatment of the function ${\rm
R}(\Delta,b,w,f)$ is the preliminary discussion on the tetrahedral
symmetries, from which the notion of log-branch emerges
straightforwardly in a natural way (see Lemma \ref{Rsym3}). Also, when
considering a $2 \leftrightarrow 3$ move configuration between
hyperbolic ideal tetrahedra, we thus have to look at all possible
branching transits, for ${\rm R}(\Delta,b,w,f)$ explicitely depends on
the branching $b$. Each branching transit leads to a specific instance of
functional lifted five term relation for ${\rm R}(\Delta,b,w,f)$, and
there actually exist some which only hold mod $(\pi^2/6)\mz$.}
\end{remark}

%% file: NSUMDINV.tex
\section{The dilogarithmic invariant $R(W,\rho)$}\label{sumdinv}

\noindent Let $\Tt=(T,b,z)$ be a $\Dd$-triangulation for $(W,\rho)$, with idealization $\Tt_{\Ii}=(T,b,w)$.
 In order to define the dilogarithmic invariant, one has to enrich $\Tt_{\Ii}$ with a system $f$ of flattenings 
{\it stable} for log-branch transits.

\begin{defi}\label{globalflat}{\rm We say that $f$ is a {\it flattening of}  $\Tt_{\Ii}=(T,b,w)$ if each $f_i$ is a 
flattening of $(\Delta_i,b_i,w_i)$ and the following global conditions are verified:
\smallskip

(1) The sum of $b$-signed log-branches around every edge of $T$ is equal to $0$.
\smallskip

(2) Condition (1) implies that $f$ verifies a similar one mod$(2)$, so that it 
represents a class $[f]\in H^1(W,\Z/2\Z)$. 
One requires that $[f]=0$.}
\end{defi}

\noindent The first condition is formally the $\log$ of the compatibility condition in Lemma \ref{compat}. 
The second ensures a nice {\it affine} structure on the set of these flattenings. A slight adaptation of a 
fundamental result due to W. Neumann \cite{N3,N2} implies the following:

\begin{teo}\label{esistenzaflatt} Any $\Ii$-triangulation $\Tt_{\Ii}=(T,b,w)$ for $(W,\rho)$ admits a flattening. 
Flattenings on $\Tt_{\Ii}=(T,b,w)$ make an affine space over a lattice. This lattice has an explicitely given basis 
made by one vector for the abstract star of each edge of $T$.
\end{teo}

\noindent Finally one can state

\begin{teo}\label{dilodinvariante} Let $\Tt_{\Ii}$ be any flattened $\Ii$-triangulation for $(W,\rho)$. Set
$$ {\rm R}(\Tt_{\Ii},f) = \sum_i *_i \ {\rm R}(\Delta_i,b_i,w_i,f_i)\in \C/(\pi^2/6)\Z\ ,$$
where $*_i=\pm 1$ according to $b_i$. The value of ${\rm R}(\Tt_{\Ii},f)$ does not depend on the choice of 
$(\Tt_{\Ii},f)$. Hence it defines an invariant ${\rm R}(W,\rho)\in \C/(\pi^2/6)\Z$ called the {\rm dilogarithmic invariant} 
of the pair $(W,\rho)$.
\end{teo}
\begin{prop}\label{CS} We have ${\rm R}(W,\rho) = {\rm CS}(\rho) + i{\rm Vol}(\rho)$ {\rm mod}$(\pi^2/6)\Z$, 
where {\rm CS} and {\rm Vol} are the Chern-Simons invariant and the volume of the character $\rho$.
\end{prop}

\noindent This proposition shows that ${\rm R}(W,\rho)$ refines the {\rm
mod}($\pi^2\Q$) dilogarithmic interpretation of the second
Cheeger-Chern-Simons class of $\rho$, due to Dupont-Sah
\cite{DS,D1}. This is in agreement (except for the {\rm
mod}$(\pi^2/6)\Z$ precision) with the results stated in \cite{N1,NY}
in the particular case when $\rho$ is the holonomy of a genuine
hyperbolic structure on $W$.

\noindent In \cite{BB3} we had roughly announced the existence (limited to $B$-characters) of 
such dilogarithmic invariant defined {\rm mod}$(\pi^2/2)\Z$. A more careful analysis of the tetrahedral 
symmetries leads to the present {\rm mod}$(\pi^2/6)\Z$ formulation.

%% file: NSUMQHI.tex
\section{Quantum dilogarithms and QHI for $(W,L,\rho)$ }\label{sumQHI}

\noindent In this section we revisit the construction of the B-QHI made in
\cite{BB1}, having as model the construction of the dilogarithmic
invariant. The result shall be a full generalization of the QHI for
arbitrary $PSL(2,\C)$-characters, and a deeper understanding of the
ultimate hyperbolic geometric nature of these invariants.
\subsection{Quantum dilogarithms}\label{Qd}

\noindent 
Let $N>1$ be a fixed odd positive integer, and put
 $\zeta=\exp(2i\pi/N)$. The \emph{Faddeev-Kashaev's quantum
 dilogarithms} were originally derived as an
 explicit matrix realization of the associator in the cyclic
 representation theory of a Borel quantum subalgebra of
 $U_{\zeta}(sl(2,\mc))$ (see \cite{FK,K1}). The matrix elements of the quantum
 dilogarithms are usually called $6j$-symbols. This derivation is
 also presented in the Appendix of \cite {BB1}, and with full details 
in \cite{B}. Here we forget this `quantum algebraic' origin, and simply
 describe and interpret the special functions one has to deal with.

\smallskip

\noindent For any complex number $x$ with $\vert x \vert <1$, consider the
analytic function $g$ defined by
$$g(x) := \prod_{j=1}^{N-1}(1 - x\zeta^j)^{j/N}$$
\noindent and set $h(x) := x^{-p}g(x)/g(1)$ when $x$ 
is non-zero (one computes that $\vert g(1) \vert = N^{1/2}$).  We
shall still denote by $g$ its analytic continuation to the complex
plane with cuts from the points $x = \exp (i\epsilon)\zeta^k,\ k
=0,\ldots,\ N-1$, $\epsilon \in \mathbb{R}$, to infinity. Hereafter we
will implicitly assume that $\epsilon$ is such that the cuts are away
from the points where $g$ is evaluated.

\noindent 
Consider the curve $\Gamma=\{x^N + y^N = z^N\} \subset \mathbb{C}P^2$
(homogeneous coordinates), and the rational functions given for any $
n \in \mathbb{N}$ by
\begin{eqnarray} \label{omeg}
\omega(x,y,z \vert n) = \prod_{j=1}^n \frac{(y/z)}{1-(x/z)\zeta^j} \quad .
\end{eqnarray}
\noindent 
These functions are periodic in their integer argument, with period
$N$. Denote by $\delta$ the $N$-periodic Kronecker symbol,
i.e. $\delta(n) = 1$ if $n \equiv 0$ mod($N$), and $\delta(n) = 0$
otherwise. Set $[x] = N^{-1}\ (1-x^N)/(1-x)$.

\medskip

\noindent 
The $N$-dimensional Faddeev-Kashaev quantum dilogarithm and its
inverse are the $N^2$-matrices whose components are the rational
functions defined on the curve $\Gamma$ by
$$R(x,y,z)_{\alpha,\beta}^{\gamma,\delta} = h(z/x)\
\zeta^{\alpha\delta+\frac{\alpha^2}{2}}\ \omega(x,y,z \vert \gamma-\alpha) \ \delta(\gamma + \delta - \beta)$$
$$\bar{R}(x,y,z)_{\gamma,\delta}^{\alpha,\beta} =  
\frac{[x/z]}{h(z/x)}\ \zeta^{-\alpha\delta-\frac{\alpha^2}{2}}\ \frac{\delta(\gamma + \delta - \beta)}{\omega(\frac{x}{\zeta},
y,z\vert \gamma-\alpha)}\ .$$

\noindent We can interpret these matrices as functions of 
$\Ii$-tetrahedra as follows. Let $(\Delta,b,w)$ be an
$\Ii$-tetrahedron. Write $w_i = -p_{i+1}/p_{i+2}$ (indices
mod($\mz/3\mz$)) as after Def. \ref {idealizable}. Recall that $p_0 +
p_1 + p_2 =0$.  Fix a common determination of the $N$-th roots of the
$p_i$'s, which we denote by $p_i'$. Set
\begin{eqnarray}
\mathfrak{L}(\Delta,b,w) = \left\lbrace \begin{array}{l}
R(p_1',p_0',-p_2') \quad {\rm if} \quad *=1 \\ \\ 
\bar{R}(p_1',p_0',-p_2') \quad {\rm if} \quad *=-1 \quad , 
\end{array} \right. \nonumber
\end{eqnarray}
\noindent 
\noindent where $*=\pm1$ according to the $b$-orientation of $\Delta$.
Note that (\ref{omeg}) implies that $\mathfrak{L}(\Delta,b,w)$ only
depends on $(b,w)$, and not on the $N$-th roots $p_i'$ of the $p_i$'s,
for it is homogeneous in these variables. One realizes that

\medskip

\noindent {\it $\mathfrak{L}(\Delta,b,w)$ does not respect the tetrahedral
symmetries.}
\medskip
 
\noindent 
Let us consider the $1$-skeleton of the cell decomposition of $\Delta$
dual to the canonical triangulation with 4 vertices. It is made by 4
edges incident to the center of $\Delta$. They are oriented by the
orientation complementary to the $b$-orientations of the dual 2-faces of
$\Delta$. Hence there is a couple of arcs incoming into $\Delta$ and a
couple of outcoming ones. One can associate to both couples a copy of
$\C^N \otimes \C^N$ (with the standard basis) and interpret
$\mathfrak{L}(\Delta,b,w)$ as a linear operator defined on the
incoming couple, with values in the outcoming one.

\noindent 
Let $\Tt_{\Ii} = (T,b,w)$ be any $\Ii$-triangulation for
$(W,\rho)$. Let us consider the $1$-skeleton $C$ of the cell
decomposition dual to $T$, with edges oriented as above. By
associating to each $(\Delta_i,b_i,w_i)$ the corresponding operator
$\mathfrak{L}(\Delta_i,b_i,w_i)$, one gets an {\it operator network}
whose complete contraction gives a scalar
$\mathfrak{L}_N(\Tt_{\Ii})\in \C$ (note that there is no edge with
free ends in $C$). This has an explicit expression as a {\it state
sum}, where the {\it states} are given by the indices of the matrices
entries as follows. A state is a function defined on the edges of $C$, with
values in $\{0,\dots, N-1\}$. For every $(\Delta_i,b_i,w_i)$, any state
$\alpha$ determines an entry (a $6j$-symbol)  
$\mathfrak{L}(\Delta_i,b_i,w_i)_\alpha$ of $\mathfrak{L}(\Delta_i,b_i,w_i)$.
Set 
$$\mathfrak{L}_N(\Tt_{\Ii})_\alpha = 
\prod_i \mathfrak{L}(\Delta_i,b_i,w_i)_\alpha $$ and
\begin{equation}\label{statesum}
\mathfrak{L}_N(\Tt_{\Ii})= 
\sum_\alpha \mathfrak{L}_N(\Tt_{\Ii})_\alpha \quad .
\end{equation}  
One realizes that

\medskip

\noindent 
{\it $\mathfrak{L}_N(\Tt_{\Ii})$ is invariant {\rm only} for some
peculiar instances of $\Ii$-transits (one being the same as for
Schaeffer's identity for the Rogers dilogarithm)}.

\medskip

\noindent 
These facts justify the following name: $\mathfrak{L}(\Delta,b,w)$ is
the $N$-dimensional {\it non symmetric quantum dilogarithm}, computed
on the given $\Ii$-tetrahedron. The special $\Ii$-transits which keep
$\mathfrak{L}_N(\Tt_{\Ii})$ invariant induce the {\it basic five term
(pentagonal) relations} satisfied by the non symmetric quantum
dilogarithm. In order to construct invariants for $(W,\rho)$ based on
the quantum dilogarithms $\mathfrak{L}$, one has to solve the same kind of
difficulties than for defining the dilogarithmic invariant, based on the classical
Rogers dilogarithm L.

\subsection{Tetrahedral symmetries}\label{symql}

\noindent 
Let $(\Delta,b,w)$ be as above.  An {\it integral charge} on
$(\Delta,b,w)$ is a $\mz$-valued map on the edges of $\Delta$ such
that $c(e)=c(e')$ for opposite edges $e$ and $e'$, and $c_0 + c_1 +
c_2 = 1$, where $c_i=c(e_i)$. Note that a charge is formally similar
to a flattening. The main difference is that the charge does not
depend on the moduli $w$. In fact a charge defines a flattening only
if $*_w = -1$. Write $N=2p+1$, and for each edge $e$ set $c'(e)=(p+1)\
c(e)$ mod($N$), viewed as a point in $\{0,\ldots,N-1\}$.
\begin{defi}\label{symqd}
{\rm The {\it symmetrized quantum dilogarithm} is the matrix-valued
operator, defined on the set of charged $\Ii$-tetrahedra
$(\Delta,b,w,c)$, given by}
\begin{eqnarray}\label{formsym}
\mathfrak{R}(\Delta,b,w,c) = \left\lbrace \begin{array}{l}

\bigl( (-p_1'/p_2')^{-c_1}\ (-p_2'/p_0')^{c_0} \bigr)^p\ R'(w \vert c) \quad {\rm if} \quad *=1 \\
\bigl( (-p_1'/p_2')^{-c_1}\ (-p_2'/p_0')^{c_0} \bigr)^p\ \bar{R}'(w \vert c) \quad {\rm if} \quad *=-1 \quad , 

\end{array} \right.
\end{eqnarray}
\noindent 
{\rm where $*=\pm1$ according to the $b$-orientation of $\Delta$, and
the matrices $R'(w \vert c)$ and $\bar{R}'(w \vert c)$ have the
components}
\begin{eqnarray}\label{symmatrix}
R'(w \vert c)_{\alpha,\beta}^{\gamma,\delta} = \zeta^{c_1'(\gamma - \alpha)} \ 
R(p_1',p_0',-p_2')_{\alpha,\beta - c_0'}^{\gamma - c_0',\delta}\nonumber \\ \\
\bar{R}'(w \vert c)_{\gamma,\delta}^{\alpha,\beta} = \zeta^{c_1'(\gamma - \alpha)}\ 
\bar{R}(p_1',p_0',-p_2')_{\gamma+c_0',\delta}^{\alpha,\beta+c_0'}\ .\nonumber
\end{eqnarray}
\end{defi}
\noindent 
Note again that (\ref{omeg}) implies that $\mathfrak{R}(\Delta,b,w,c)$
only depend on $(b,w,c)$, and not on the choice of the $N$-th roots $p_i'$ of the
$p_i$'s. Here we identify $1/2 \in \mz/N\mz$ with $p+1$.

\begin{remark}\label{newsym}{\rm 
This symmetrization sligtly differs from the one adopted in
\cite{BB1} for $B$-characters, see Remark \ref{Bsym} below. 
This new one has been obtained by
formally copying the flattening's contributions in $ \exp({\rm
R}(\Delta,b,w,f))$, postulating that charges and flattenings should
essentially play the same role. In fact, it works.}
\end{remark} 

\noindent 
Recall that the permutation group on four elements, which is the
symmetry group of a branched abstract tetrahedron, is generated by the
transpositions $(01)$, $(12)$ and $(23)$.  The following lemma
describes the tetrahedral symmetries of $\mathfrak{R}$.

\begin{lem} \label{6jsym} 
Let $(\Delta,b,w,c)$ be a charged $\Ii$-tetrahedron.  Changing the
orientation of the edge $e_0$, $e_1$ or $e_0'$ by the transposition
$(01)$, $(12)$ or $(23)$ of its vertices gives respectively
\begin{eqnarray}
\mathfrak{R}\bigl((01)(\Delta,b,w,c)\bigr) & \equiv_{(\pm \zeta^{\Z})} & T_1^{-1}\ \mathfrak{R}(\Delta,b,w,c) \ T_1 
\nonumber \\
\mathfrak{R}\bigl((12)(\Delta,b,w,c)\bigr) & \equiv_{(\pm \zeta^{\Z})} &  S_1^{-1}\ \mathfrak{R}(\Delta,b,w,c)\ T_2 
\nonumber \\
\mathfrak{R}\bigl((23)(\Delta,b,w,c)\bigr) & \equiv_{(\pm \zeta^{\Z})} &  S_2^{-1}\ \mathfrak{R}(\Delta,b,w,c)\ S_2\quad   
\nonumber ,
\end{eqnarray}
\noindent  where $\equiv_{(\pm \zeta^{\Z})}$ means equality up to sign and multiplication by $N$-th roots of unity. Here we write $T_1=T \otimes 1$, etc..., and $T$ and $S$ are the $N$-dimensional invertible square 
matrices with components $T_{m,n} = \nu\ \zeta^{\frac{m^2}{2}}\delta(m + n)$ and $S_{m,n} = N^{-\frac{1}{2}}\zeta^{mn}$, where 
$\nu = g(1)/ \vert g(1) \vert$.
\end{lem}

\noindent The proof is based on the same computations that gave Prop. 9.6 in the Appendix of \cite{BB1}, adapted to the new symmetrization.

\subsection{Complete five term relations}

\noindent Replacing the non symmetric quantum dilogarithms with the symmetrized ones in (\ref{statesum}) for {\it charged} $\Ii$-triangulations
$(T,b,w,c)$ of $(W,\rho)$, one obtains state sums
$\mathfrak{R}_N(\Tt_{\Ii},c)$. The next step is to complete the
$\Ii$-transits to suitable {\it charged} $\Ii$-transits, in order to
realize the full transit invariance of the state sum's value.  As the
charges are moduli-independent, also their transit is.  

\noindent One says that there is a {\it charge transit} $(T,c)
\leftrightarrow (T',c')$ if $c'$
equals $c$ on the edges of the abstract tetrahedra of $T$ not involved
in the move, and for any other edge $e$ we have the {\it transit of
sum} condition:
$$\sum_{a \in \epsilon_T^{-1}(e)} c(a) = \sum_{a' \in
\epsilon_{T'}^{-1}(e)} c'(a')\quad .$$ Note that, for $2 \to 3$
transits, this implies that the sum of the charges around the new edge
after the move is equal to $2$.  Note also that a charge transit
coincides with a flattening transit, providing that the signs $*_b$'s
and $*_w$'s satisfy certain conditions which are easy to determine.
\noindent 
By using Lemma \ref{6jsym} and the basic five term (pentagonal)
relations for the non symmetric quantum dilogarithms we finally get:
\begin{prop}\label{6jtransit} 
For any charged $2 \leftrightarrow 3$ $\Ii$-transit $(T,w,c) \leftrightarrow
(T',w',c')$ we have
$$\prod_{\Delta_i \subset T} \ \mathfrak{R}(\Delta_i,b,w,c)
\equiv_{(\pm \zeta^{\Z})}\ \prod_{\Delta_i' \subset T'}
\ \mathfrak{R}(\Delta_i',b',w',c')\quad .$$
\end{prop}

\subsection{Link-fixing and QHI for triples $(W,L,\rho)$}

\noindent As for flattenings in the definition of the dilogarithmic
invariant, the next step is to give the right notion of global
integral charge which should be stable for charged $\Ii$-transits. The
solution of this problem is much more elaborated than for the
dilogarithmic invariant. The naive idea should be to require that the
sum of charges around each edge of $T$ is equal to $2$. But simple
combinatorial considerations show that such tentative global integral
charges do not exist.  A way to overcome this difficulty is to fix an
arbitrary link $L$ in $W$ (considered up to ambient isotopy) and to
inglobe it in all the construction of the QHI. 

\noindent A $\Dd$-triangulation
$(\Tt,H)$ for $(W,L,\rho)$ is a $\Dd$-triangulation $\Tt$ for
$(W,\rho)$ with a Hamiltonian (i.e. which contains all the vertices of
$T$) subcomplex $H$ which realizes the link $L$. An integral charge
$c$ on $(\Tt,H)$ satisfies by definition the following conditions:
\smallskip

(1) the sum of its values equal $2$ (resp. $0$) around the edges of $T
    \setminus H$ (resp. $H$);
\smallskip

(2) it satisfies the same $\Z/2\Z$-cohomological condition than for flattenings: there is a class $[c] \in H^1(W;\mz/2\mz)$, and we require that $[c]=0$. 
\smallskip

\noindent 
Note that any charge $c$ eventually encodes $H$, hence the link $L$.
The existence of charged $\Dd$- and $\Ii$-triangulations $(\Tt,c)$,
$(\Tt_{\Ii},c)$ for $(W,L,\rho)$ is a rather demanding fact proved in
\cite{B} and \cite{BB1}. Again, the good affine structure on the set of charges is a variation of the fundamental Neumann's result Th. \ref{esistenzaflatt}.
The presence of the link must be integrated also in the transits;
in particular the
condition to be Hamiltonian must be preserved. For instance, for $2
\rightarrow 3$ transits the new edge must be in $T' \setminus H'$; for
positive bubble moves one edge of the initial triangle must belong to
$H$, and it is replaced with the two edges connected to it and to the
new vertex in $T'$. Finally one achieves the construction of the QHI:

\begin{teo}\label{QHIinv} The value of the (normalized) state sum 
$$H_N(T_{\Ii},c) = N^{-n_0}\ \mathfrak{R}_N(\Tt_{\Ii},c)\quad ,$$ 
where $n_0$ is the number of vertices of $T$, does not
depend on the choice of the charged $\Ii$-triangulation $(T_{\Ii},c)$
for $(W,L,\rho)$, up to sign and multiplication by $N$-th roots of
unity.  Hence, up to this ambiguity, for any triple $(W,L,\rho)$ it
defines a \emph{quantum hyperbolic invariant} $H_N(W,L,\rho) \in \mc$.
\end{teo}
\begin{remark}\label{linkprob} {\rm The presence of the link $L$ and the ambiguity up to
multiplication by $N$-th roots of unity of $H_N(W,L,\rho)$ are
entirely a consequence of the specific symmetrization procedure of the
quantum dilogarithms we have adopted. We guess that suitable
variations of this procedure (by using moduli-dependent charges, or by
following even more closely the strategy used for uniformizing the
Rogers dilogarithm) allows one to define the QHI directly for
$(W,\rho)$, and also to well define them only up to multiplication by
$\exp(i\pi/12)$. We postpone this problem to a future work.}
\end{remark} 

\begin{remark} \label{Bsym} {\rm The QHI for $B$-characters defined
in \cite{BB1} and their state sum formulas differ from those in
\ref{QHIinv} by a scalar factor depending on the cocycle $z$ of $\Tt$,
not only on $\Tt_{\Ii}$. This is a consequence of the different
symmetrization adopted in \cite{BB1}. There, it consisted in replacing in
(\ref{formsym}) the scalar state-independent factor in front of the matrices $R'$ and $\bar{R}'$ by $(-q_2')^p$,
where the $q_j$'s have been defined in Remark
\ref{idealrem} (3) and $'$ denotes, as before, a common determination
of the $N$-th roots of the $q_j$'s. Let us denote by $\mathfrak{R}^B(\Tt,c)$ the associated state sum. 

\noindent Then, the statement of Lemma \ref{6jsym} is unchanged, except that there is no sign ambiguity. 
However, in Prop. \ref{6jtransit} one has to multiply both sides by the respective
$\textstyle Q_2 := \prod_i (-q_2')_i^p$. It is a remarkable but somewhat fortuitous fact that, 
for $B$-characters and for any positive $2 \to 3$ $\Dd$-transit $\Tt \rightarrow \Tt'$, one has 
$Q_2(\Tt')/Q_2(\Tt)=x(e)^{2p}$, where $x(e)$ is the upper-diagonal value of the cocycle $z$ on the 
new edge in $T' \setminus H'$. Normalizing $\mathfrak{R}^B(\Tt,c)$ by dividing it 
by $\textstyle \prod_{e \in T \setminus H} x(e)^{2p}$, it eventually gives a well defined invariant up 
to $N$-th roots of unity. The same procedure for general $PSL(2,\mc)$-characters (using the $p_2'$'s 
instead of the $q_2'$'s) does not seem to work, due to the fact that the explicit formula for $P_2(\Tt')/P_2(\Tt)$ 
heavily depends on the branching.}
\end{remark} 

%% file: SUMSCISSORS.tex
\section{Scissors congruence classes}\label{sumscissors}

\noindent Consider the free $\Z$-modules $\Z(\Dd)$ and $\Z(\Ii)$ which are
respectively generated by all $\Dd$-tetrahedra and all $\Ii$-tetrahedra.  Let
$\Pp(\Dd)$ (resp. $\Pp(\Ii)$) be the quotient of $\Z(\Dd)$
(resp. $\Z(\Ii)$) by all instances of the relations associated to the
tetrahedral symmetries in Lemma \ref {Iisym}, and to the
$2\leftrightarrow 3$ $\Dd$-transits (resp. $\Ii$-transits). 

\noindent We call $\Pp(\Dd)$ (resp. $\Pp(\Ii)$) the $\Dd$- (resp. $\Ii$-) {\it (pre)-Bloch group}.  
Prop. \ref{DdomI} implies that there is a surjective
homomorphism $I: \Pp(\Dd) \to \Pp(\Ii)$.

\noindent Working with
either flattened or charged $\Ii$-tetrahedra, and using the
appropriate lifts of the tetrahedral symmetries and of the five term
identities, one obtains the $F$- and $C$-{\it refined}
$\Ii$-(pre)-Bloch groups $\Pp(\Ii)_F$ and $\Pp(\Ii)_C$.
\smallskip

\noindent 
Clearly, every $\Dd$-triangulation $\Tt$ and every $\Ii$-triangulation
$\Tt_{\Ii}$ for $(W,\rho)$ defines an element $\cG(\Tt) \in
\Pp(\Dd)$ and $\cG(\Tt_{\Ii})\in \Pp(\Ii)$ respectively. Similarly,
any flattened $( \Tt_{\Ii},f)$ defines an element $\cG(\Tt_{\Ii},f)\in
\Pp(\Ii)_F$, and any charged $( \Tt_{\Ii},c)$ for $(W,L,\rho)$ defines an element $\cG(\Tt_{\Ii},c)\in \Pp(\Ii)_C$. 
One has

\begin{prop}\label{classisciss}  The elements $\cG(\Tt)$, $\cG(\Tt_{\Ii})$, 
$\cG(\Tt_{\Ii},f)$ and $\cG(\Tt_{\Ii},c)$ do not depend on the choices in their respective arguments. 
Hence they define invariants $\cG_D(W,\rho)\in  \Pp(\Dd)$, $\cG_I(W,\rho)\in  \Pp(\Ii)$, $\cG_F(W,\rho)\in  \Pp(\Ii)_F$ 
and
$\cG_C(W,L,\rho)\in  \Pp(\Ii)_C$. 

\noindent These are generally called {\rm scissors congruence classes} of $(W,\rho)$ or of $(W,L,\rho)$. 
Moreover, the homomorphism $I$ sends $\cG_D(W,\rho)$ onto $\cG_I(W,\rho)$.
\end{prop}

\noindent Clearly, the formula in Th. \ref{dilodinvariante} defines
a function $ {\rm R}: \Pp(\Ii)_F \to \C/(\pi^2/6)\Z$ such that
$$ {\rm R}(W,\rho) = {\rm R}(\cG_F(W,\rho))\ .$$ 
\noindent One would like to interpret also
the QHI as evaluations on $\cG_C(W,L,\rho)$ of suitable functions
defined on $\Pp(\Ii)_C$. This is roughly true, but there is some
subtilities (see the discussion at the end of Section 5 of \cite{BB3}).

\noindent Lemma \ref {compat} implies that every  $\Ii$-triangulation $\Tt_{\Ii}$ for $(W,\rho)$
 also defines a class $\alpha(\Tt_{\Ii}) \in
H_3(PSL(2,\C),\Z)$, where $PSL(2,\C)$ is endowed with the discrete
topology. One can prove that also $\alpha(\Tt_{\Ii})$ does not depend
on the choice of $\Tt_{\Ii}$, hence it defines an invariant $\alpha
(W,\rho) \in H_3(PSL(2,\C),\Z)$. This is a main ingredient of the
group-cohomological approach of Dupont-Sah.
 

%% file: SUMCOMMENT.tex
\section{Further comments on the proofs}\label{commenti}
The proofs of Th. \ref{dilodinvariante} and Th. \ref{QHIinv} have
the very same structure. In fact one follows, almost {\it verbatim}, the proof given in 
$\{$Sections 1-6 $+$ the Appendix$\}$ of \cite{BB1} for the $B$-QHI. 
For the dilogarithmic invariant the proof is easier, because one does not have
to keep track of the link $L$ during the transits.

\noindent A key point is that quasi-regular triangulations can be connected by
quasi-regular transits. This is proved in Prop. 2.10 of \cite{BB1}. This implies that also $\Dd$-transits 
and $\Ii$-transits are {\it generically} possible, and this is enough.  In fact, another delicate point is the invariance w.r.t. 
the flattenings or the charges. This
heavily depends on Neumann's crucial Th. \ref{esistenzaflatt} and runs like in the proof of Th. 4.6 in  \cite{BB1}. 

\noindent For scissors congruence classes, the above {\it generic} existence
of $\Dd$-transits and $\Ii$-transits is not enough. We have to refine Prop. 2.10 of \cite{BB1} (by using the same technics indeed), in order to get
the full existence of transits. There is another subtility: the
(pre)-Bloch groups are defined by only using the relations coming from
$2 \leftrightarrow 3$ transits. But to connect triangulations (with
possibly a different number of vertices), one has to use also the
`bubble' transits. So one has to check that this induces relations
which are in fact consequences of the other ones. This is proved as in
Corollary 4.5 of \cite {BB1}.

\noindent The situation is a little more complicated if one works with arbitrary branchings, 
or if one allows also {\it degenerate} ideal tetrahedra. General branchings were already treated in 
\cite {BB1}, so one simply repeats those arguments. 
In presence of degenerate tetrahedra, the invariance of the dilogarithmic sum 
w.r.t. the branching becomes more delicate. This is 
automatic under the non-degenerate assumption, thanks to Lemma \ref {Rsym3}. 
Anyway, one can obtain it by using the quasi-regular refinement of
the following result of F. Costantino \cite{C}:
\medskip

\noindent {\it Two branchings on a same given triangulation $T$ of $W$ can be connected by transits
of branchings.}
\medskip

\noindent This allows to reduce also the branching invariance to the invariance by transits.

%% file: NSUMVOL.tex
\section{A Volume Conjecture for the QHI of $(W,L,\rho)$}\label{sumvol}
It is clear from the discussion in Section \ref{sumQHI} that the structure
of the QHI is modeled on the one of $\exp((1/2i\pi){\rm R}(W,\rho))$.  
In fact one introduces the $1/2i\pi$-factor to have a perfect
agreement of the behaviour with respect to the involution $W \to -W$, 
given by the change of orientation of $W$:
\begin{lem}\label{-W} 
Denote by $^*$ is the complex conjugation. We have
$$\begin{array}{c}
H_N(-W,L,\rho) = \bigl(H_N(W,L,\rho^*)\bigr)^* \\
\exp((1/2i\pi){\rm R}(-W,\rho))=\bigl(\exp((1/2i\pi){\rm
R}(W,\rho^*))\bigr)^* \quad .
\end{array}$$
\end{lem}
The first identity is proved like in Prop. 6.1. in \cite{BB1}. The
second one is a consequence of Prop. \ref{CS}, and of ${\rm
CS}(\rho^*)= {\rm CS}(\rho)$ and ${\rm Vol}(\rho^*)= -{\rm
Vol}(\rho)$. These structural coincidences and the actual asymptotic
behaviour of the quantum dilogarithms (see below) motivate the
following {\it Volume Conjecture} for the asymptotic
expansion of QHI, when $N\to \infty$:
\begin{conj}\label{volconj} 
There exist invariants $C=C(W,L,\rho) \in \mc \ {\rm mod}(\pi^2/6)\mz$ and $D=D(W,L,\rho)\in \mc^*$ such that
$$ \bigl( H_N(W,L,\rho)\bigr)^N \equiv \left[ \exp\left( \frac{C+ N {\rm R}(W,\rho)}{2i\pi}\right) \right]^{N} \left(D +
 \mathcal{O}(\frac{1}{N})\right)$$  
\noindent where $\equiv$ means equality up to multiplication by the integer
powers of $\exp(i\pi/12)$.
\end{conj}
\noindent Conjecture \ref{volconj} says at first that
$H_N(W,L,\rho)^N$ has an exponential growth rate. Assuming it, the
fact that $\exp (C/i\pi)$, $\exp({\rm R}(W,\rho)/i\pi)$ and $D$ are
well-determined invariants of $(W,L,\rho)$ follows from the invariance
of $H_N(W,L,\rho)^N$ and the uniqueness of the coefficients of
asymptotic power series expansions.

\noindent At present, the nature of
$C$ and $D$ is somewhat mysterious to us. There are no reasons to
expect that, for instance, $C=0$ or $D=1$. We have expressed the conjecture in terms of the $N$-th power of $H_N(W,L,\rho)$ so as to kill its multiplicative ambiguity up to $N$-th roots of unity. (The statement is formally the same as in \cite{BB3}, except that it was given mod$(i\pi/4)\mz$, see the discussion at the end of Section \ref{sumdinv}). 

\smallskip

\noindent Classical manipulations of one-variable complex analysis with the so-called
Faddeev's {\it non-compact dilogarithm} \cite{F} allows one to prove
that when $N \rightarrow \infty$ (we use the functions introduced in section \ref{Qd}):
\begin{equation}\label{versdilo}
g(z/x)\ \omega(x,y,z \vert n) \sim  (y/z)^n\ \ \exp \biggl[\frac{N}{2i\pi}\ \bigl({\rm Li}_2((x/z)\zeta^n)) +\log(x/z)^2 -\pi\log(x/z)  +\pi^2 \bigr) \biggr]\nonumber\ ,
\end{equation}
where log is, as before, the standard branch of the logarithm. Rewriting $n$ in terms of states and charges, ones derives from this formula the leading term of the asymptotic expansion of the symmetrized quantum dilogarithms.
This corroborates Conjecture \ref{volconj}. 
\begin{remarks}\label{nolink}{\rm 
The above Volume Conjecture predicts, in particular, that the
dominant term of the asymptotic expansion of the QHI for $N \to
\infty$ is {\it not}-sensitive to the link-fixing used to get one
specific global symmetrization of the state sums. This is not clear
for $C$ and $D$. There is a simple way to modify the dilogarithmic
invariant so as to make it link-sensitive. For that, it is enough to
set
$$ {\rm R}(W,L,\rho) = {\rm R}(W,\rho) + (i\pi/2)\log({\rm Tr}(\alpha (L))\quad ,$$ 
where $\alpha$ is any representative of $\rho$ and $L$ is arbitrarily oriented, and considered as an element of the fundamental group of
$W$. Note that ${\rm R}(W,L,\rho)$ does not depend on the choices we
made. It can be computed by using any $\Dd$-triangulation
$\Tt$ for $(W,L,\rho)$, but not only by using the idealization $\Tt_{\Ii}$.}
\end{remarks}

%% file: SUMCUSPED.tex
\section{Cusped manifolds}\label{sumcusped}

\noindent This section is less definitive than the rest of the paper. The final achievement of the results 
presented below is stricly related to the solution of the problem mentionned in Remark \ref{linkprob}.

\subsection{Dilogarithmic invariant and QHI}
Let $M$ be an oriented complete non-compact and finite volume
hyperbolic 3-manifold (shortly: $M$ is a {\it cusped
manifold}). Fix a triangulation of $M$ made by
embedded geodesic ideal tetrahedra. For simplicity, assume that it is branched. 
This is not really necessary here, due to the symmetry relations in Lemma \ref{Rsym3} and Lemma \ref{6jsym}; 
on the contrary, for the case of $(W,\rho)$ treated before, the branching was part of the structure of 
$\Dd$-triangulations, which dominate the $\Ii$-triangulations.  

\noindent With the usual notations, this triangulation of $M$ can
be represented as $\Tt_{\Ii}=(T,b,w)$ such that every
$(\Delta_i,b_i,w_i)$ is {\it quasi-geometric}, i.e. geometric in the sense of Def. \ref{geometrico} or 
possibly degenerate.

\noindent One can enrich $\Tt_{\Ii}$ with suitable
flattenings $f$ or integral charges $c$ (again this is essentially due
to \cite{N3}). Their definition is similar to the one given for
$(W,\rho)$, but there is a further condition which mimics the usual
{\it completeness condition} satisfied by the moduli. So one can
define both ${\rm R}(\Tt_{\Ii},f) \in \C/(\pi^2/6)\Z$ and
$H_N(\Tt_{\Ii},c)$.

\noindent One would like to prove that they do not depend on any choice,
 hence that they define invariants for $M$, just by using our direct methods
 based on $\Ii$-transits. It is not clear to us if any two such enriched quasi-geometric $\Ii$-triangulations 
of $M$ can be connected by enriched $\Ii$-transits. However, for ${\rm R}(\Tt_{\Ii},f)$ one can prove it 
if we impose for instance the following further restriction on $\Tt_{\Ii}$:
 
\begin{prop}\label{invdcusp} Suppose that $M$ admits triangulations $(T,b)$ such that the volume 
function, defined on the set of all $\Ii$-triangulations supported by $(T,b)$, has a unique maximum at $\Tt_{\Ii}$. 
Then ${\rm R}(\Tt_{\Ii},f)$ does not depend on any choice, thus defining a dilogarithmic invariant 
${\rm R}(M) \in \mc$ {\rm mod}$(\pi^2/6)\mz$.
\end{prop}

\noindent Here, the set of all $\Ii$-triangulations supported by $(T,b)$ 
includes {\it all} (non necessarily quasi-geometric) solutions of the
compatibility and completeness equations supported by $T$. It is not
known to us if any cusped $M$ has such special triangulations, but, as
is used in the Snappea program it often happens. 

\noindent The same result should hold true for the expected link-free definition of QHI 
(see Remark \ref{linkprob}); then, let us assume that also the QHI $H_N(M)$ of $M$ are well-defined. 
The Volume Conjecture \ref{volconj} for $(W,\rho)$ can be repeated {\it verbatim} for these invariants of
cusped manifolds.

\subsection{About the Volume Conjecture for the colored Jones invariants
$J_N(L)$} Recall that Kashaev's Volume Conjecture \cite{K2}, reformulated
in terms of the (suitably normalized) colored Jones
invariants $J_N(L)$ of links in $S^3$ by Murakami-Murakami, states that for an hyperbolic
knot $L$ we have
$$ \lim_{N\to \infty} (2\pi/N) \log ( |J_N(L)|) = {\rm Vol}(M) \quad ,$$
where $M$ is the cusped manifold given by the hyperbolic complement of
$L$ in $S^3$. Due to Remark \ref{nolink}, one has:
\smallskip

\noindent 
{\it The Volume Conjecture \ref{volconj} for $H_N(W,L,\rho)$
is not compatible with the one stated above for $J_N(L)$ if
one also assumes (as it is currently done, see {\rm \cite{K3}}) that $J_N(L)^N$ 
coincides with $H_N(S^3,L,\rho_0)^N$, where $\rho_0$
is the necessarily trivial character.}
\medskip

\noindent We guess that: 

\smallskip

\noindent (1) the last assumption above is not correct (we will elaborate
on this point in \cite{BB2}).

\smallskip

\noindent (2) The Volume Conjecture for $J_N(L)$ makes sense only if it could be related to Conjecture \ref{volconj} for general cusped
manifolds. At least they are formally compatible. However, for the
moment, we do not see a systematic way to identify (even
asymptotically) $J_N(L)$ and $H_N(M)$, where $M$ is the cusped
manifold $S^3\setminus L$.

\smallskip

\noindent (3) A consistent relationship between all these Volume Conjectures
could be obtained thanks to Thurston's {\it hyperbolic Dehn filling} and
{\it double limit procedure}. Indeed, let $M$ be a cusped manifold and
$(W_n,L_n,\rho_n)$ be a sequence of compact hyperbolic Dehn fillings
of $M$ converging to $M$. Here, $L_n$ denotes the link made
of the short simple geodesics in $W_n$ forming the {\it cores} of the fillings, and $\rho_n$ is
the holonomy of the hyperbolic manifold $W_n$. Presumably, one has
$$ {\rm R}(W_n,\rho_n) \longrightarrow {\rm R}(M)\quad ,\ n\to \infty \ , $$
and for every fixed $N$ also
$$ H_N(W_n,L_n,\rho_n) \longrightarrow H_N(M)\quad ,\ n\to \infty \ .$$
So, following Conjecture \ref{volconj}, by taking a `double limit' we are led to:
 
\begin{conj}\label{2volconj} The dominant term of the asymptotic expansion when $n,\ N \to \infty$ of $\bigl( H_N(W_n,L_n,\rho_n)\bigr)^N$  is equal to
$\exp\bigl((N^2/2i\pi){\rm R}(M)\bigr)$, up to multiplication by
integer powers of $\exp(i\pi/12)$.
\end{conj}

\noindent  The Volume Conjecture for $J_N$ seems to be corroborated
by few numerical computations, obtained via a quite formal and
`optimistic' use of the {\it stationary phase method}. We also guess
that the above considerations could help to have a correct understanding
of the {\it meaning} of those formal manipulations.